\documentclass[11pt]{article}
\usepackage{tikz}
\usepackage{amsmath}
\usepackage{amsfonts}
\usepackage{amssymb}
\usepackage{mathtools}
\usepackage{graphicx}
\usepackage{geometry}

\usetikzlibrary{arrows.meta, decorations.markings, calc, positioning}
\newcommand{\R}{\mathbb{R}}
\begin{document}
\begin{center}
\large\bf 
Contraction Mapping: Case Studies
\\
\vspace{0.2cm}
\end{center}
\begin{center}
\small{Shamanth Sreekanth}\\
\small{Electrical and Systems Engineering, University of Pennsylvania}
\end{center}

\vspace{0.2cm}

\section*{Abstract}

While exploring dynamical systems, we often come across the principle of contraction mapping, or better known as the Banach fixed point theorem. It is an essential concept based on successive approximations, whose utility comes from two main guarantees: establishing existence and uniqueness of a solution, and establishing constructive proof. The intent of this manuscript is to break down two major proofs incorporating this in ordinary differential equations (ODEs), and make them a little more understandable step-by-step to an audience that presumably has adequate knowledge of modern calculus and real analysis. These are not original proofs, only original narration.

\section{The Picard-Lindelöf Theorem}

Here, we address a fundamental question in ODEs. Suppose we have a system that evolves in time, and we are given only the following conditions:
\begin{itemize}
\item a well-defined starting point: \(y(x_0) = y_0\)
\item how the function changes near that point: \( \frac{dy}{dx} = f(x, y)\)
\end{itemize}
They form what is called an initial value problem (IVP). It defines how the system evolves and anchors the system’s trajectory with a fixed point, which, as a result, serves as an important prerequisite for prediction and control. Now, how do we know for sure if 
\begin{itemize}
\item a solution really exists for those conditions, and
\item the solution is unique (i.e., there aren't multiple conflicting ways the system can evolve in)?
\end{itemize} The Picard-Lindelöf theorem says, if \(f(x, y)\) is
\begin{itemize}
\item  continuous around the point \( (x_0, y_0) \), and
\item Lipschitz continuous with respect to \( y \), 
\end{itemize} then,
\begin{itemize}
\item there exists one solution \( y(x) \) that solves the differential equation near \( x_0 \), and
\item this solution is unique.
\end{itemize}

\noindent\textbf{Definition 1} \textit{Continuity} \\
\textit{Continuity ensures that the function behaves predictably and smoothly in a neighborhood around the point, without jumps or breaks.
A function \( f : D \subseteq \mathbb{R}^2 \to \mathbb{R} \)  is said to be continuous at a point \( (x_0, y_0) \in D \) if, as the point \( (x, y) \) approaches \( (x_0, y_0) \), the value of \( f(x, y) \) approaches \( f(x_0, y_0) \). Basically, for every \( \varepsilon > 0 \), there exists a \( \delta > 0 \) such that
\[
\sqrt{(x - x_0)^2 + (y - y_0)^2} < \delta \implies
|f(x, y) - f(x_0, y_0)| < \varepsilon
\]}
\noindent\textbf{Definition 2} \textit{Lipschitz Continuity} \\
\textit{This implies that the rate of change of \( f \) with respect to \( y \) is uniformly bounded across the domain. Lipschitz continuity is stronger than ordinary continuity and is essential for guaranteeing uniqueness of solutions to differential equations. Basically, fixing a point on the x-axis, we vary \(y\), to see how sharply the function changes.
A function \( f : D \subseteq \mathbb{R}^2 \to \mathbb{R} \) is said to be Lipschitz continuous in the variable \( y \) if there exists a constant \( L > 0 \) such that for all \( (x, y_1), (x, y_2) \in D \),}
\[
|f(x, y_1) - f(x, y_2)| \leq L |y_1 - y_2|
\]
With this in hand, we can now state the theorem in full mathematical rigor.
\subsection*{\bf Theorem:}
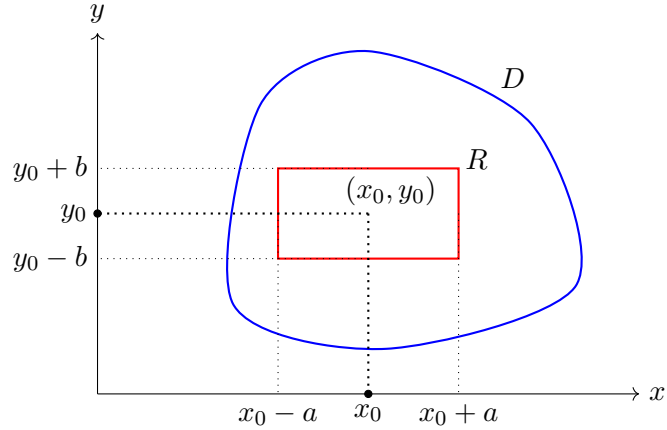
\begin{center}
\begin{tikzpicture}[scale=1.2]
  \draw[->] (0,0) -- (6,0) node[right] {$x$};
  \draw[->] (0,0) -- (0,4) node[above] {$y$};

  \def\xo{3}
  \def\yo{2}
  \def\a{1}
  \def\b{0.5}

  \draw[red, thick] ({\xo - \a},{\yo - \b}) rectangle ({\xo + \a},{\yo + \b});
  \node at ({\xo + \a + 0.2}, {\yo + \b + 0.1}) {$R$};

  \draw[blue, thick] plot [smooth cycle] coordinates {
    (1.5,1) 
    (1.8,3.2) 
    (3,3.8) 
    (4.8,3) 
    (5.3,1.2) 
    (3.2,0.5)
  };
  \node at (4.6, 3.5) {$D$};

  \draw[dotted, thick] (\xo,\yo) -- (\xo,0); 
  \draw[dotted, thick] (\xo,\yo) -- (0,\yo); 

  \draw[fill=black] (\xo,0) circle (0.04); 
  \node[below] at (\xo, 0) {$x_0$};

  \draw[fill=black] (0,\yo) circle (0.04); 
  \node[left] at (0, \yo) {$y_0$};

  \draw[dotted] ({\xo - \a}, \yo) -- ({\xo - \a}, 0);
  \node[below] at ({\xo - \a}, 0) {$x_0 - a$};

  \draw[dotted] ({\xo + \a}, \yo) -- ({\xo + \a}, 0);
  \node[below] at ({\xo + \a}, 0) {$x_0 + a$};

  \draw[dotted] (\xo, {\yo - \b}) -- (0, {\yo - \b});
  \node[left] at (0, {\yo - \b}) {$y_0 - b$};

  \draw[dotted] (\xo, {\yo + \b}) -- (0, {\yo + \b});
  \node[left] at (0, {\yo + \b}) {$y_0 + b$};

  \node at ({\xo + 0.25}, {\yo + 0.25}) {$(x_0, y_0)$};

\end{tikzpicture}
\end{center}

Let \( D \) be a domain in \( \mathbb{R}^2 \), and let \( f : D \to \mathbb{R} \)
be a real-valued function that satisfies the following conditions:
\begin{itemize}
    \item[(i)] \( f \) is continuous on \( D \),
    \item[(ii)] \( f(x, y) \) is Lipschitz continuous with respect to \( y \) on \( D \), with Lipschitz constant \( L \geq 0 \).
\end{itemize}

Let \( (x_0, y_0) \) be an interior point of \( D \).
Let \( a > 0 \), \( b > 0 \) be some constants such that the rectangle
\[
R = \left\{ (x, y) \in \mathbb{R}^2 \mid |x - x_0| \leq a,\ |y - y_0| \leq b \right\} \subset D
\]\\
Define
\[
M = \max_{(x, y) \in R} f(x, y), \quad h = \min\left(a, \frac{b}{M}\right)
\]

\noindent
Then, the initial value problem
\[
\frac{dy}{dx} = f(x, y), \quad y(x_0) = y_0
\]
\noindent
has a unique solution \( y(x) \) on the interval \( |x - x_0| \leq h \).\\

\subsection*{Proof:}
Let's look at the statement a little carefully. \(D \subseteq \mathbb{R}^2\) is an open set in a 2-dimensional plane where our function \(f(x,y)\) exists. It's essentially our playground for now. The more important thing to focus on is \(R\), which a rectangle such that its breadth is 2\(a\) and its height is 2\(b\), with the center of the rectangle being (\(x_0, y_0\)). We assume it houses the two prerequisite properties of the theorem: \(f(x,y)\) being continuous and Lipschitz continuous, thereby bounding our function to it. This becomes our "workspace".

We define \(h\) similarly. It is the maximum interval \( |x - x_0| \leq h \) inside which the solution is guaranteed to exist. Therefore it is essential for \(h\) to be defined in a way where the function stays within \(R\) and doesn't go out of bounds. We obviously know \(a\) is the most we can move horizontally in either direction without leaving \(R\). But the vertical change isn't directly set. It depends on how fast the solution can possible move, i.e., \(M = \max f(x, y)\).
\[
\frac{dy}{dx}\leq M
\]
\[
dy\leq Mdx
\]
\[
|y-y_0|\leq M|x-x_0|
\]
\[
|y-y_0|\leq Mh
\]
But, \(|y-y_0|\leq b\) vertically. Because of this hard bound, we choose \(h\) as something small enough such that 
\[
Mh \leq b
\]
\[
h \leq \frac{b}{M}
\]
So essentially, \(h\) is either \(a\), if the function allows for it, or is \(\frac{b}{M}\), whichever is minimum.
\[\quad h = \min\left(a, \frac{b}{M}\right)\]

Okay, so now we understand what we're working on, and what's on the table. We know a starting value. And we know how the function behaves. We know the solution could theoretically exist in that neighborhood. But we don't know the solution itself, or how it's going to look like, or if it even exists. Therefore, we take the initial condition and we try moving along the direction dictated by \(f(x,y)\) to hopefully shape ourselves a solution \(y(x)\). \\ Let \(x = x_0  +h\) for some \(h >0\).
\[
y'=\frac{dy}{dx}=f(x,y)
\]
\[
\int_{x_0}^{x}y(t)dt = \int_{x_0}^{x}f(t, y(t))dt
\]
\[
y(x) - y(x_0) = \int_{x_0}^{x}f(t, y(t))dt
\]
\[
y(x) = y(x_0) + \int_{x_0}^{x}f(t, y(t))dt
\]
This equation is the integral form of the IVP, and it helps us construct a recursive sequence which we'll follow to approximate ourselves to the final solution \(y(x)\). We'll go ahead and denote every iteration of this sequence as \(\phi\). \\Our starting value will be
\[
\phi_0(x) = y_0
\]
Plugging this into consecutive iterations,
\begin{align*}
\phi_1(x) &= y_0 + \int_{x_0}^{x} f(t, \phi_0(t))\, dt \\
\phi_2(x) &= y_0 + \int_{x_0}^{x} f(t, \phi_1(t))\, dt \\
\vdotswithin{=} \\
\phi_n(x) &= y_0 + \int_{x_0}^{x} f(t, \phi_{n-1}(t))\, dt \tag{1} \label{piterate}
\end{align*}

This sequence (\ref{piterate}) defined for [\(x_0, x_0+h\)] is called Picard's Iterates, and it forms the heart of our proof. To validate ourselves that we're using the right setup, we'll first prove that the set {\(\phi_n\)} is well defined and continuous, and subsequently that \(f(x, \phi_n(x))\) is also well defined and continuous such that \(|\phi_n(x) - y_0| \leq b\). Basically we prove that all values of {\(\phi_n\)} lie within \(R\) \\

\noindent\textbf{Definition 3} \textit{Well-defined \\
A function or operation is said to be well-defined if it assigns a unique output to every valid input in its domain. If \( f : A \to B \), and \( a \in A \), \( b \in B \) then the value \( f(a) = b\) is unique to \(a\). In this context, (\ref{piterate}) is said to be well-defined if the function \( f(t, \phi_{n-1}(t)) \) is integrable on [\(x_0, x_0+h\)], and the resulting integral yields a valid, unique real number.}\\

We can use mathematical induction to do this, as it's a neat way to establish that a particular statement holds true for all \(n \in \mathbb{N}\) iterations of a sequence. It consists of two steps - a \textit{base case} where we show the statement holds true for a value like \(n = 0\) or \(n = 1\), and an \textit{inductive case} where we assume it holds true for some arbitrary \(n=k\) (a hypothesis), and use it to prove it holds true for \(n=k+1\).

Assume \(\phi_n(x)\) exists, has continuous derivative on [\(x_0, x_0+h\)], and \(|\phi_n(x) - y_0| \leq b\). This naturally implies \((x, \phi_n(x)) \in R\). Then, \(f(x, \phi_n(x))\) is well-defined and continuous, as \(f(x, \phi_n(x)) \leq M\).
Now if,
\[
\phi_{n+1}(x) = y_0 + \int_{x_0}^{x} f(t, \phi_{n}(t)) \, dt
\]
\[
\phi_{n+1}(x) - y_0 = \int_{x_0}^{x} f(t, \phi_{n}(t)) \, dt
\]
\[
|\phi_{n+1}(x) - y_0| = \left| \int_{x_0}^{x} f(t, \phi_{n}(t)) \, dt \right|
\leq \int_{x_0}^{x} \left| f(t, \phi_{n-1}(t)) \right| dt\]

\[
\leq \int_{x_0}^{x} M \, dt = M(x - x_0)
\]

\[
\leq Mh \quad \text{since } x - x_0 \leq h
\leq b
\]

Because of this, \((x, \phi_n(x))\) lies in \(R_1\), so \(f(x, \phi_{n+1}(x))\) is defined and continuous on \([x_0, x_0 + h]\)\\

When \(n = 1\)\\
\[\phi_1(x) = y_0 + \int_{x_0}^{x} f(t, \phi_0(t)) \, dt
\]

Since \(\phi_0(x) = y_0\) and \(y_0\) are well-defined and continuous, \(\phi_1(x)\) is also well-defined \& continuous.

\[
|\phi_1(x) - y_0| \leq \int_{x_0}^{x} |f(t, y_0)| \, dt \leq M(x - x_0)
\leq Mh \leq b
\]

Thus, by method of induction, the sequence (\ref{piterate}) possesses all desired properties in \([x_0, x_0 + h]\).

What we have here is a series of approximations of what could be our solution, becoming consecutively more accurate with every successive iteration. Let's compare two successive iterates. 
\[
|\phi_{n+1}(x) - \phi_n(x)| = \left| \int_{x_0}^{x} f(t, \phi_{n}(t)) \, dt - \int_{x_0}^{x} f(t, \phi_{n-1}(t)) \, dt \right|
\leq \int_{x_0}^{x} |f(t, \phi_{n}(t))- f(t, \phi_{n-1}(t)) |\, dt
\]
Notice how the integrand terms are just the function \(f(x,y)\) evaluated at two different points on the y-axis, keeping the x-axis value constant. This is precisely what Lipschitz continuity is, as defined previously. If we take L as the Lipschitz constant,
\[
|\phi_{n+1}(x) - \phi_n(x)|\leq \int_{x_0}^{x} L|\phi_{n}(t)- \phi_{n-1}(t) |\, dt
\]
Let's now compute the starting few iterations and see what we get.\\
Taking \(n=0\):
\[
|\phi_{1}(x) - \phi_0(x)|\leq \int_{x_0}^{x} |f(t,\phi_0(t)|\, dt\leq\int_{x_0}^{x}Mdt=M(x-x_0) \tag{2} \label{0iterate}
\]
Taking \(n=1\):
\[
|\phi_{2}(x) - \phi_1(x)|\leq \int_{x_0}^{x} L|\phi_{1}(t)- \phi_{0}(t) |\, dt
\]
From (\ref{0iterate}),
\[
\int_{x_0}^{x} L|\phi_{1}(t)- \phi_{0}(t) |\, dt \leq\int_{x_0}^{x} LM(t-x_0)\, dt 
\]
\[
=LM\int_{x_0}^{x}(t-x_0)dt
\]
\[
=LM\left[\frac{(t-x_0)^2}{2}\right]_{x_0}^{x}
\]
\[
=LM\frac{(x-x_0)^2}{2} \tag{3} \label{1iterate}
\]
Taking \(n=2\):
\[
|\phi_{3}(x) - \phi_2(x)|\leq \int_{x_0}^{x} L|\phi_{2}(t)- \phi_{1}(t) |\, dt
\]
From (\ref{1iterate}),
\[
\int_{x_0}^{x} L|\phi_{2}(t)- \phi_{1}(t) |\, dt \leq
\int_{x_0}^{x}L\left(LM\frac{(t-x_0)^2}{2}\right)\, dt
\]
\[
=L^2M\int_{x_0}^{x} \frac{(t-x_0)^2}{2}\, dt
\]
\[
=L^2M\left[\frac{(t-x_0)^3}{6}\right]_{x_0}^x
\]
\[
=L^2M\frac{(x-x_0)^3}{6}
\]
We see a pattern emerging. As \(n \xrightarrow{}\infty\),
\[
|\phi_{n+1}(x) - \phi_n(x)|\leq L^{n-1}M\frac{(x-x_0)^n}{n!}
\]
We know from previous definitions that \(x-x_0\leq h\). Rearranging the terms a bit,
\[
\leq L^nL^{-1}M\frac{h^n}{n!}
\]
\[
=\frac{L^n}{L}\frac{Mh^n}{n!}
\]
\[
=\frac{M}{L}\frac{(Lh)^n}{n!} \tag{4} \label{finalform}
\]

Now what's left is for us to prove that the series converges such that the elements of \(\phi_n(x)\) get closer to each other, to give us a final approximation. This is called Cauchy convergence, and is the modus operandi of the Picard's iterates.\\

\noindent\textbf{Definition 4} \textit{Cauchy Convergence} \\
A sequence \( \{\phi_n\} \) in a metric space is said to be \textit{Cauchy convergent} (or simply, a \textit{Cauchy sequence}) if for every \( \varepsilon > 0 \), there exists an number \( N \in \mathbb{N} \) such that for all \( m, n \geq N \), we have
\[
|x_n - x_m| < \varepsilon.
\]
This basically means, beyond a certain iteration count \(N\), any two iterations of that sequence (they don't even have to be consecutive), their difference will be less than \(\varepsilon\). In this context, the sequence of functions \( \{\phi_n(x)\} \) is Cauchy convergent in the space of continuous functions on \([x_0, x_0 + h]\) in the supremum norm, if
\[
\sup_{x \in [x_0, x_0 + h]} |\phi_n(x) - \phi_m(x)| < \varepsilon \quad \text{for all } m,n \geq N.
\]
Although the difference is very subtle, a good question to ask is why supremum? Why not maximum?\\

\noindent\textbf{Definition 5} \textit{Supremum Norm} \\
Given a function \( f \) defined on a closed interval \([a, b]\), the \textit{supremum norm} (also called the \textit{uniform norm}) of \( f \) is defined as
\[
\|f\|_{\infty} = \sup_{x \in [a, b]} |f(x)|.
\]
It represents the maximum absolute value attained by the function on the interval \([a, b]\). For a sequence of functions \( \{\phi_n(x)\} \), convergence in the supremum norm means that
\[
\|\phi_n - \phi\|_{\infty} = \sup_{x \in [a, b]} |\phi_n(x) - \phi(x)| \to 0 \quad \text{as } n \to \infty.
\]

The supremum (sup) and maximum (max) are very similar, but the key difference lies in whether the maximum value is actually attained. Maximum is a value that a function actually reaches on the interval. Supremum is the smallest possible value that we can classify as an upper bound, even if the function never actually hits that value. What we basically mean to say is, the minimum possible upper bound of \(|\phi_n(x) - \phi_m(x)|\) is 0, even if it never truly reaches 0 and gets very close to it. Because the series has this nature, this possibility, of never truly reaching 0, we use supremum, as a way of generalization.\\

So, by definition of Cauchy convergence, let's take two iterations beyond a point \( N \in \mathbb{N} \). Assuming \(m > n\), we can break them down to their consecutive parts:
\[
|\phi_m(x)-\phi_n(x)| \leq |\phi_m(x)-\phi_{m-1}(x)| + |\phi_{m-1}(x)-\phi_{m-2}(x)| + ... + |\phi_{n+1}(x)-\phi_n(x)|
\]
We know \ref{finalform} holds true for consecutive iterates. Applying it as a summation for the above equation,
\[
|\phi_m(x)-\phi_n(x)| \leq \frac{M}{L}\sum_{k=n}^{m-1}\frac{(Lh)^k}{k!}
\tag{5}\label{powertail}\]
We make a key observation here, something that solves one half of our problem statement. (\ref{powertail}) is the tail of a well-known power series of the form \(e^z\). If \(z=Lh\),
\[e^z = \sum_{k=0}^{\infty}\frac{z^k}{k!}\]
\[e^{Lh} = \sum_{k=0}^{\infty}\frac{Lh^k}{k!}\]
We know the function \(e^{Lh}\) converges because it is well known to do so. It can be verified through a ratio test, if needed.
There exists a small non-zero value \(\varepsilon>0\), and \( N \in \mathbb{N} \) such that, for any of \(N_0\geq N\), 
\[\frac{M}{L}\sum_{k=N_0}^{\infty}\frac{Lh^k}{k!} < \varepsilon\]
As \(n\xrightarrow{}\infty\), \(\phi_n(x)\) gets increasingly closer to being our hypothesized solution \(\phi(x)\). To evaluate this, we take their difference.
\[
|\phi_n(x) -\phi(x)| = \left|\int_{x_0}^xf(t, \phi_{n-1}(t))dt - \int_{x_0}^xf(t, y(t))dt\right|
\]
\[
\leq\int_{x_0}^x L|(\phi_{n-1}(t)-\phi(t))|dt
\]
We've already established that \(|\phi_{n-1}-\phi| \leq ||\phi_{n-1} - \phi||\), and that the supremum \(||\phi_{n-1} - \phi||\) is a constant upper bound value.
\[
\leq L\int_{x_0}^{x}||\phi_{n-1}-\phi||dt = L||\phi_{n-1}-\phi||\int_{x_0}^{x}dt
\]
\[
=L||\phi_{n-1}-\phi||(|x-x_0|)
\]
\[
<Lh||\phi_{n-1}-\phi|| \rightarrow 0
\]
When \(n\) gets sufficiently large the difference between \(\phi_n\) and \(\phi\) converges further, making \(||\phi_{n-1}-\phi|| \rightarrow0\) within some \(\varepsilon>0\).

Hence we have proved that there indeed exists a solution.

To check uniqueness, an easy way is to assume there are multiple solutions, and work our way backwards.
Suppose there are two solutions, \(y_1(x)\) and \(y_2(x)\). They both have the same IVP where,
\[
y_1(x_0) = y_0
\]
\[
y_2(x_0) = y_0
\]
\[
y_1(x)=y_0+\int_{x_0}^{x}f(t,y_1(t))dt
\]
\[
y_2(x)=y_0+\int_{x_0}^{x}f(t,y_2(t))dt
\]
Subtracting them from each other,
\[
|y_1(x)-y_2(x)| = \left|\int_{x_0}^{x}f(t,y_1(t))dt-f(t,y_2(t))dt\right|
\]
\[
\leq\int_{x_0}^x L|(y_{1}(t)-y_2(t))|dt
\]
If
\[
z(x):=|y_1(x)-y_2(x)|
\]
Then
\[
z(x)\leq L\int_{x_0}^{x}z(t)dt
\]
We observe that this form is similar to the Grönwall's Inequality definition.\\

\noindent\textbf{Lemma 1} \textit{(Grönwall's Inequality)} \\
Let \(z:[a,b] \to \mathbb{R}\) be a continuous, nonnegative function. Suppose there exists a constant \(C \geq 0\) and a continuous, nonnegative function \(A:[a,b] \to \mathbb{R}\) such that
\[
z(x) \leq C + \int_a^x A(t)\,z(t)\,dt, \quad \forall x \in [a,b].
\]
Then \(z(x)\) is bounded by
\[
z(x) \leq C \exp\!\Bigg(\int_a^x A(t)\,dt\Bigg), \quad \forall x \in [a,b].
\]

\noindent This lemma basically tells that if a function is bounded by a constant and an integral term that grows proportionally to itself, then the function can grow at most exponentially, with the rate determined by the function \(A(t)\).
In our case however, \(A(t)=L\) and \(C=0\), which results in
\[
z(x)\leq 0(e^{L\int_{x_0}^{x}dt})=0(e^{L(x-x_0)})=0
\]
Because \(|y_1(x)-y_2(x)| = z(x)=0\), this implies that there \textbf{exists} only one \textbf{unique} solution \(y(x)\), which finally completes the proof. \\

By proving the existence and uniqueness of a solution, this theorem guarantees that for a given state, there is only one specific future. Without this guarantee, physical prediction would be impossible, as systems could arbitrarily follow different paths from the same starting point.

\newpage
\section{The Hartman Grobman Theorem}

Having established via Picard-Lindelöf that a unique solution for an ODE exists, we now turn to a theorem that approximates its behavior. As explicit analytic solutions for nonlinear systems are rarely obtainable, we typically rely on linearization, approximating the system near a fixed point. The Hartman-Grobman Theorem asserts that near a hyperbolic fixed point, the nonlinear flow is topologically equivalent to its linearized counterpart. This guarantees the existence of a homeomorphism $H$ that continuously deforms the curved nonlinear trajectories into the predictable paths of the linear system, to determine local stability.

\subsection*{Theorem:}

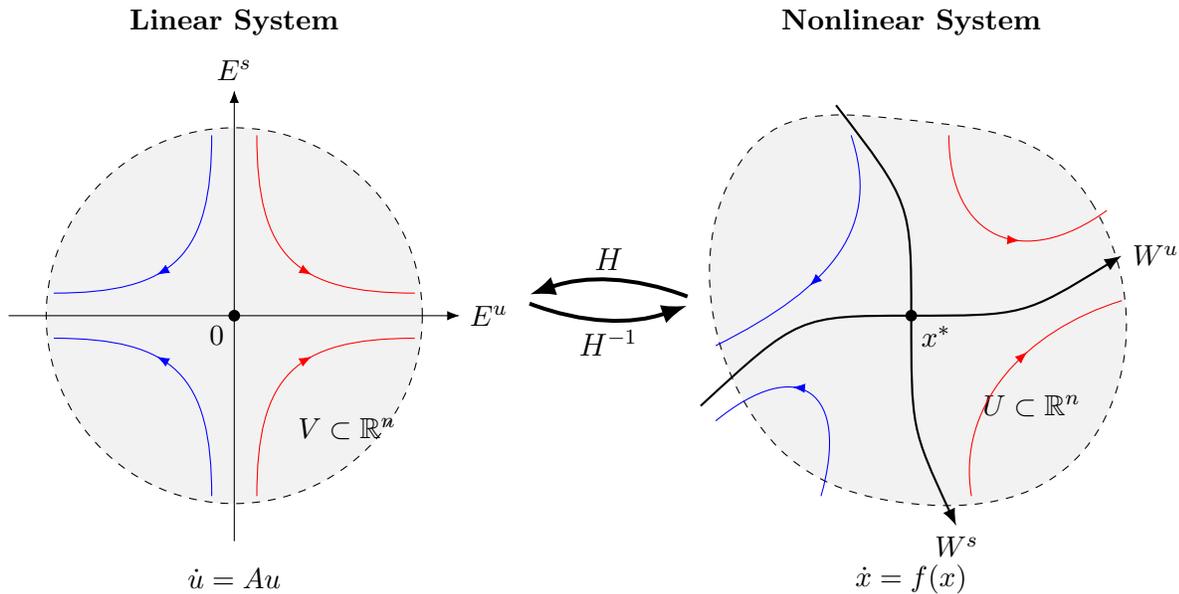
\begin{tikzpicture}[
    >=Latex,
    flow/.style={
        decoration={
            markings,
            mark=at position 0.6 with {\arrow{>}}
        },
        postaction={decorate}
    }
]

\def\offset{9}

\begin{scope}[local bounding box=LinearScope]
    \draw[dashed, fill=gray!10] (0,0) circle (2.5);
    \node at (1.5, -1.5) {$V \subset \mathbb{R}^n$};
    
    \draw[->] (-3,0) -- (3,0) node[right] {$E^u$};
    \draw[->] (0,-3) -- (0,3) node[above] {$E^s$};
    
    \draw[flow, blue] (-0.3, 2.4) .. controls (-0.3, 0.8) and (-0.8, 0.3) .. (-2.4, 0.3);
    \draw[flow, red] (0.3, 2.4) .. controls (0.3, 0.8) and (0.8, 0.3) .. (2.4, 0.3);
    \draw[flow, blue] (-0.3, -2.4) .. controls (-0.3, -0.8) and (-0.8, -0.3) .. (-2.4, -0.3);
    \draw[flow, red] (0.3, -2.4) .. controls (0.3, -0.8) and (0.8, -0.3) .. (2.4, -0.3);
    
    \filldraw (0,0) circle (2pt) node[below left] {${0}$};
    
    \node[above, font=\bfseries] at (0, 3.6) {Linear System};
    \node at (0, -3.5) {$\dot{u} = Au$};
\end{scope}

\begin{scope}[shift={(\offset,0)}, local bounding box=NonlinearScope]
    \draw[dashed, fill=gray!10] plot [smooth cycle, tension=0.8] coordinates {
        (0, 2.6) (2.4, 1.5) (2.5, -1.5) (0, -2.5) (-2.3, -1.0) (-2.3, 2.0)
    };
    \node at (1.6, -1.2) {$U \subset \mathbb{R}^n$};

    \draw[->, thick] (-2.8, -1.2) .. controls (-1.5, 0) .. (0,0) .. controls (1.5, 0) .. (2.8, 0.8) node[right] {$W^u$};
    
    \draw[->, thick] (-1.0, 2.8) .. controls (0, 1.5) .. (0,0) .. controls (0, -1.5) .. (0.6, -2.8) node[below] {$W^s$};
    
    
    \draw[flow, red] (0.5, 2.4) .. controls (0.5, 1.0) and (1.5, 0.6) .. (2.6, 1.4);
    
    \draw[flow, blue] (-0.8, 2.4) .. controls (-0.4, 1.2) and (-1.0, 0.4) .. (-2.6, -0.4);
    
    \draw[flow, blue] (-1.2, -2.4) .. controls (-0.8, -1.0) and (-1.5, -0.5) .. (-2.6, -1.4);
    
    \draw[flow, red] (0.8, -2.4) .. controls (0.6, -1.2) and (1.5, -0.2) .. (2.8, 0.2);

    \filldraw (0,0) circle (2pt) node[below right] {${x}^*$};

    \node[above, font=\bfseries] at (0, 3.6) {Nonlinear System};
    \node at (0, -3.5) {$\dot{x} = f(x)$};
\end{scope}

\draw[->, line width=1.5pt, shorten >= 5pt, shorten <= 5pt] 
    (NonlinearScope.west) to[bend right=20] node[midway, above, font=\large\bfseries] {$H$} 
    (LinearScope.east);

\draw[->, line width=1.5pt, shorten >= 5pt, shorten <= 5pt] 
    (LinearScope.east) to[bend right=20] node[midway, below, font=\large\bfseries] {$H^{-1}$} 
    (NonlinearScope.west);

\end{tikzpicture}

Let ${x} \in \R^n$. \\Consider a non-linear system $\dot{{x}} = f({x})$ with flow $\phi_t$.
Consider a linear system $\dot{{x}} = A{x}$, where $A = Df({x}^*)$ (Jacobian matrix of $f({x})$ at ${x}^*$). Assume ${x}^*$ is a hyperbolic fixed point and ${x}^* = {0}$. Let $f$ be $C^1$ (differentiable at $f$, continuous at $f$) on some $E \subset \R^n$ that contains ${0}$. Let $I_0 \subset \R$, $U \subset \R^n$, $V \subset \R^n$ with ${0} \in I_0, U, V$.

Then, there exists a homeomorphism $H: U \to V$ such that, $\forall$ initial points ${x}_0 \in U \subset E \subset \R^n$, and $\forall t \in I_0$, $ H \circ \phi_t({x}_0) = e^{At} H({x}_0) $\\

\noindent\textbf{Definition 6} \textit{Flow}\\
A \textbf{flow}, denoted $\phi(t; {x}_0)$, is the function that describes the evolution of a dynamical system over time. If $\dot{{x}} = F({x})$ is the system of differential equations, the flow $\phi(t; {x}_0)$ gives the state (position) of the system at time $t$ given that it started at the initial condition ${x}_0$.\\

\noindent\textbf{Definition 7} \textit{Homeomorphism}\\
A function $h: X \to Y$ is said to be a \textbf{homeomorphism} if it is a continuous bijection (meaning it is one-to-one and onto) and its inverse function, $h^{-1}$, is also continuous. The existence of a homeomorphism means the spaces $X$ and $Y$ are topologically equivalent.\\

\noindent\textbf{Definition 8} \textit{Hyperbolic Fixed Point}\\
A fixed point of a dynamical system is said to be \textbf{hyperbolic} if all the eigenvalues of the system's Jacobian matrix, when evaluated at that fixed point, have non-zero real parts.\\

\noindent\textbf{Definition 9} \textit{Jacobian}\\
The \textbf{Jacobian} matrix (or "total derivative") is a matrix containing all the first-order partial derivatives of a vector-valued function. If $f: \mathbb{R}^n \to \mathbb{R}^m$, its Jacobian $J$ is an $m \times n$ matrix where the entry in row $i$, column $j$ is $\frac{\partial f_i}{\partial x_j}$.

\subsection*{Proof}

Essentially, we are trying to find a magic box ``Homeomorphism $H$'' that continuously deforms the complicated curved trajectories of the nonlinear system into the simple straight line (or spiral) trajectories of the system.

Assume ${x}^* = 0$.

The A matrix is put in a block diagonal form (also called Jordan form (\textit{Perko, 1991})).
$$ A = \begin{pmatrix} P & 0 \\ 0 & Q \end{pmatrix} $$

where, \\
$P$ = matrix of stable dynamics (sink). $\lambda_P = a+ib$, $a < 0$ \\
$Q$ = matrix of unstable dynamics (source). $\lambda_Q = a+ib$, $a > 0$ \\

Consider state vector ${x}$, it is also split into
$$ {x}_0 = \begin{pmatrix} {y}_0 \\ {z}_0 \end{pmatrix} $$
$$ \phi_t({x}_0) = {x}(t, {x}_0) = \begin{pmatrix} {y}(t, {y}_0, {z}_0) \\ {z}(t, {y}_0, {z}_0) \end{pmatrix} $$
where $\phi_t({x}_0)$ is the true nonlinear flow of soln. at time $t$. \\

Here, ${y}_0 \in E^s$ is the stable component, belonging to the stable subspace $E^s$ spanned by the eigenvectors of $P$. ${z}_0 \in E^u$ is the unstable component belonging to the unstable subspace $E^u$ spanned by the eigenvectors of $Q$.

Now, instead of solving the continuous time problem, let us look at one discrete time step $t=1$.
Define two functions $\tilde{Y}$ and $\tilde{Z}$ such that,
\begin{equation}
    \begin{aligned}
    \tilde{Y}({y}_0, {z}_0) &= {y}(1, {y}_0, {z}_0) - e^{P}{y}_0\\
    \tilde{Z}({y}_0, {z}_0) &= {z}(1, {y}_0, {z}_0) - e^{Q}{z}_0 
    \end{aligned} 
    \tag{6}
\end{equation}
$\tilde{Y}$ and $\tilde{Z}$ basically measure the error terms between the true nonlinear flow $\phi_t({x}_0)$ and linear systems (hypothesized).
$e^P {y}_0$ is the predicted position of the stable component after $t=1$, had it followed the linear system $\dot{{y}} = P{y}$, whose solution is ${y}(t) = e^{Pt}{y}_0$.
Vice versa for $\tilde{Z}$.
e know that if ${z}_0 = {0}$, then ${y}_0 = {z}_0 = {0}$.
$$ \tilde{Y}({0}) = \tilde{Z}({0}) = {0} $$
At the fixed point ${0}$, the linear approximation is perfect, so the error is ${0}$.
This implies that $D\tilde{Y}({0}) = D\tilde{Z}({0}) = 0$ as well.
Since $f$ is $C^1$ on $E$, by (1) and (2), $\tilde{Y}$ and $\tilde{Z}$ are also $C^1$ on $E$.
Because of this property, we say, that for a certain $a \in \R$, there exists some radius $S_0$ around the origin, small enough that the norm of the derivatives of the error terms $D\tilde{Y}$ and $D\tilde{Z}$ inside that radius are bounded by $a$.
\begin{equation}
    \begin{aligned}
        \| D \tilde{Y}(y_0, z_0) \| &\leq a \\
        \| D \tilde{Z}(y_0, z_0) \| &\leq a
    \end{aligned}
    \tag{7}
\end{equation}
To find a global solution, we will use successive approximations. However, $\tilde{Y}$ and $\tilde{Z}$ cannot be used as they are undefined outside $S_0$. We create new smooth functions $Y$ and $Z$, such that if $|{y}|^2 + |{z}|^2 > S_0^2$ then $Y = Z = 0$, and if $|{y}_0|^2 + |{z}_0|^2 \le (S_0/2)^2$ then $Y = \tilde{Y}$ and $Z = \tilde{Z}$.

This ensures that as ${x}$ moves outward from the origin towards the boundary of $S_0$, we get a smooth transition. It is globally well behaved and locally identical to the real problem.

If $\sqrt{|{y}_0|^2 + |{z}_0|^2}$ is the L2 norm, then by Mean-Value Theorem, we get,
\begin{align*}
    |Y| &\le a \sqrt{|{y}_0|^2 + |{z}_0|^2} \\
    \sqrt{|{y}_0|^2 + |{z}_0|^2} &< |{y}_0| + |{z}_0| \quad \text{(by triangle inequality)} \\
    |Y| &\le a \sqrt{|{y}_0|^2 + |{z}_0|^2} \le a(|{y}_0| + |{z}_0|)
\end{align*}
Likewise,
$$ |Z| \le a \sqrt{|{y}_0|^2 + |{z}_0|^2} \le a(|{y}_0| + |{z}_0|) $$
Here, $a$ can be considered the Lipschitz constant, when $|Y|$ or $|Z|$ is measured from the origin.

We know, that $e^{P}$ represents the linear flow of the stable subspace, and $e^Q$ represents the linear flow of the unstable subspace.
$$ B = e^P \qquad C = e^Q $$

Define $b$ and $c$ in such a way that,
\begin{equation}
    \begin{aligned}
    b &= \|B\| < 1 \\
    c &= \|C^{-1}\| < 1
    \end{aligned}
    \tag{8}
\end{equation}

Given that the flows traverse back to the origin (fixed point).
We will use all of the above to start a contraction mapping to prove the homeomorphism $H: U \to V$; $H \circ T = L \circ H$ exists.

$L$, $T$ and $H$ are defined as follows -
$$ L({y}, {z}) = \begin{bmatrix} B{y} \\ C{z} \end{bmatrix} $$
This is the linear map of the flow, hypothesized.

$$ T({y}, {z}) = \begin{bmatrix} B{y} + \tilde{Y}({y}, {z}) \\ C{z} + \tilde{Z}({y}, {z}) \end{bmatrix} $$
This is the full complete non linear flow consisting of the linear map + the error term.

$$ H({y}, {z}) = \begin{pmatrix} \phi({y}, {z}) \\ \psi({y}, {z}) \end{pmatrix} $$
Basically, $L({x}) = e^{At}{x}$ and $T({x}) = \phi_t({x})$.

Then, $H \circ T = L \circ H$ is equivalent to -
\begin{equation}
    \begin{aligned}   
    \phi(B{y} + Y({y}, {z}), C{z} + Z({y}, {z})) &= B\phi({y}, {z})\\
    \psi(B{y} + Y({y}, {z}), C{z} + Z({y}, {z})) &= C\psi({y}, {z}) 
    \end{aligned}
    \tag{9}
\end{equation}
This is analogous to letting $T$ play out and then applying $H$ to it, as compared to applying $H$ first and then following through with $L$.

Let us choose the second equation of (9) for our successive approximation. We do this as it naturally allows us to form an iterative contraction.
\begin{equation}
    \psi({y}, {z}) = C^{-1} \psi(B{y} + Y({y}, {z}), C{z} + Z({y}, {z})) \tag{10}
\end{equation}
The goal is to arrive at a function $\psi$ that satisfies (10).
Let us generate a sequence of functions, where each function is a better approximation than the last.
$$ \psi_0({y}, {z}) = {z} $$
${z}$ represents the unstable part of ${x}$.

\noindent Likewise from (10),
\begin{equation*}
    \psi_{k+1}(y,z) = C^{-1} \psi_k(B y + Y(y,z), Cz + Z(y,z))
\end{equation*}

\noindent By induction, it follows that $\psi_k$ is continuous because $\phi_t$ is continuous. Outside $S_0$, for some arbitrary radius (say $2s_0$), $\psi_k(y,z) = z$ for $|y| + |z| \ge 2S_0$. This is because the error term equivalent $(Y, Z)$ is defined as $Y, Z = 0$ for $|x| > S_0$.

\bigskip

\noindent The proof from \textit{Perko (1991)} bounds two successive iterations by
\begin{equation*}
    |\psi_j(y,z) - \psi_{j-1}(y,z)| \le M r^j (|y| + |z|)^\delta
\end{equation*}
We can arrive at this by means of induction.

\noindent For $j=1$,
\begin{align*}
    |\Psi_1({y}, {z}) - \Psi_0({y}, {z})| &= |C^{-1}\Psi_0(B{y} + Y, C{z} + Z) - {z}| \\
        &= |C^{-1}(C{z} + Z({y}, {z})) - {z}| \\
    &= |( {z} + C^{-1}Z({y}, {z}) ) - {z}| \\
    &= |C^{-1}Z({y}, {z})| \\
    &\le ||C^{-1}|| \cdot |Z({y}, {z})|\\
    &\le c \cdot a(|{y}| + |{z}|) \tag{Using the bounds $||C^{-1}|| \le c$ and $|Z| \le a(|{y}| + |{z}|)$}
\end{align*}
We know that in the region of interest, $(|{y}| + |{z}|) \le 2s_0$, which implies $(|{y}| + |{z}|)^{1-\delta} \le (2s_0)^{1-\delta}$ since $1-\delta > 0$.
\begin{align*}
    ca(|{y}| + |{z}|)^1 &= ca(|{y}| + |{z}|)^{\delta} (|{y}| + |{z}|)^{1-\delta} \\
    &\le ca \cdot (|{y}| + |{z}|)^{\delta} \cdot (2s_0)^{1-\delta} \\
    &= \left( \frac{ac(2s_0)^{1-\delta}}{r} \right) \cdot r \cdot (|{y}| + |{z}|)^{\delta} \\
    &= M r (|{y}| + |{z}|)^{\delta} \tag{11}
\end{align*}
Thus, the inequality holds for $j=1$.

\noindent The term '$r$' represents overall steepness of the function by taking $a$, $b$ and $c$ from previous equations to calculate a `worst-case' bound.
\begin{equation*}
    r = c \left[ 2 \max(a, b, c) \right]^\delta \tag{12}
\end{equation*}
$c$ is a shrinking function ($c = ||C^{-1}||$). The `worst-case' steepness is represented by $[2 \max(a, b, c)]$ by mixing the linear parts $b$ and $c$, and the non-linear error $a$, which consists of $Y$ and $Z$, therefore represented by $2a$. This term is easier to use in the rest of the proof. $\delta$ is a helper function to ensure Hölder continuity in (11) and convergence guarantee in (12). $M$ is basically just a starting constant to stabilize the bound to make the inequality true for $j=1$. \\

\noindent Now assume the inequality holds for $j=k$:
$$ |\Psi_k({y}, {z}) - \Psi_{k-1}({y}, {z})| \le M r^k (|{y}| + |{z}|)^\delta $$
We must show it holds for $j=k+1$. 
\begin{align*}
    |\Psi_{k+1}({y}, {z}) - \Psi_k({y}, {z})| &= |C^{-1}\Psi_k(B{y} + Y, C{z} + Z) - C^{-1}\Psi_{k-1}(B{y} + Y, C{z} + Z)| \\
    &\le ||C^{-1}|| \cdot |\Psi_k(B{y} + Y, C{z} + Z) - \Psi_{k-1}(B{y} + Y, C{z} + Z)| \\
    &\le c \cdot |\Psi_k(B{y} + Y, C{z} + Z) - \Psi_{k-1}(B{y} + Y, C{z} + Z)|
\end{align*}
\begin{align*}
    &\le c \cdot \left( M r^k (|B{y} + Y({y}, {z})| + |C{z} + Z({y}, {z})|)^\delta \right) \tag{Substituting $({y}', {z}') = (B{y} + Y, C{z} + Z)$} \\
    &\le c M r^k \left( b|{y}| + a(|{y}| + |{z}|) + c|{z}| + a(|{y}| + |{z}|) \right)^\delta \tag{Using bounds for B, C, Y, Z} \\
    &= c M r^k \left( (b+a)|{y}| + (c+a)|{z}| + a(|{y}| + |{z}|) \right)^\delta \\
    &\le c M r^k \left( 2\max(a, b, c) (|{y}| + |{z}|) \right)^\delta \\
    &= M r^k \cdot {c [2\max(a, b, c)]^\delta} \cdot (|{y}| + |{z}|)^\delta \\
    &= M r^{k+1} (|{y}| + |{z}|)^\delta
\end{align*}

This completes the proof by induction. The proposition is true for all $j \ge 1$, and therefore we have hence constructed the unstable component of the homeomorphism. The final stage of the proof involves constructing the stable component of the homeomorphism, denoted as $\Phi({y}, {z})$, which requires an adjustment regarding the direction of time. For the unstable component $\Psi$, we solved the equation $C\Psi = \Psi(T)$ by rearranging it to $\Psi = C^{-1}\Psi(T)$, as $C$ is an expanding function, so its inverse $C^{-1}$ is a contraction ($||C^{-1}|| < 1$). For the stable component $\Phi$, the governing equation is \(B\Phi({y}, {z}) = \Phi(T({y}, {z}))\). Since $B$ represents stable dynamics (contraction), its inverse $B^{-1}$ represents expansion ($||B^{-1}|| > 1$). Iterating with an expanding operator would cause the approximations to diverge rather than converge. To resolve this, we run time backwards. We utilize the inverse of the nonlinear map, $T^{-1}$, which exists provided the nonlinear error (Lipschitz constant $a$) is sufficiently small. The inverse map is defined as:
\[
T^{-1}({y}, {z}) = \begin{bmatrix} B^{-1}{y} + {Y}_1({y}, {z}) \\ C^{-1}{z} + {Z}_1({y}, {z}) \end{bmatrix}
\]
Here, ${Y}_1$ and ${Z}_1$ represent the nonlinear error terms for the backward time step. By rewriting the original functional equation using the inverse map $T^{-1}$, we obtain:
\[
B^{-1}\Phi({y}, {z}) = \Phi(T^{-1}({y}, {z}))
\]
\[
\Phi({y}, {z}) = B \Phi(T^{-1}({y}, {z}))
\]
This transformation works well because the operator on the outside is now $B$. Since $B$ corresponds to the stable subspace, we know that $b = ||B|| < 1$. Thus, the operator is a contraction, guaranteeing convergence. We apply the method of successive approximations using the following recursive sequence:
\begin{align*}
    \text{Initial Guess:} \quad & \Phi_0({y}, {z}) = {y} \\
    \text{Iteration:} \quad & \Phi_{k+1}({y}, {z}) = B \Phi_k(T^{-1}({y}, {z}))
\end{align*}
\begin{equation*}
    |\Phi_j(y,z) - \Phi_{j-1}(y,z)| \le M r^j (|y| + |z|)^\delta
\end{equation*}
Because $b < 1$, this sequence converges uniformly to the unique solution $\Phi({y}, {z})$. Having solved for both the stable component $\Phi$ and the unstable component $\Psi$, we construct the final homeomorphism $H$ by stacking the two functions:
\[
H({y}, {z}) = \begin{bmatrix} \Phi({y}, {z}) \\ \Psi({y}, {z}) \end{bmatrix}
\]
This function $H$ serves as the topological conjugacy between the nonlinear system and its linearization, completing the proof of the Hartman-Grobman Theorem.

\section{Conclusions}

It is striking that two different goals, anchoring a single trajectory (Picard) and linearizing a global flow (Hartman) boil down to the same Banach Fixed Point mechanism. Whether the task is shrinking the integral error in Picard’s iterates or bounding the error terms between flows as \textit{Hartman (1960)} \cite{hartman1960} originally demonstrated, the iterative mechanics are identical.

It is worth mentioning that Hartman followed a different approach in his original proof, by first establishing linearization for a discrete map using successive approximations, then extending this to continuous flows using an integral averaging technique, which necessitated a stronger $C^2$ assumption. In contrast, this breakdown (following \textit{Perko (1991)} \cite{perko1991}) streamlines this by immediately decomposing the system into stable ($P$) and unstable ($Q$) blocks. Instead of averaging flows, we solve for the stable and unstable components of the homeomorphism separately, using the forward map $T$ for the unstable part and the inverse map $T^{-1}$ for the stable part, to ensure contraction while requiring only $C^1$ smoothness. One particular novelty this manuscript achieves is the explanation of how the the upper bound for the successive iterates ($|\psi_j(y,z) - \psi_{j-1}(y,z)| \le M r^j (|y| + |z|)^\delta$) came into existence, which is skipped in Perko's textbook.

The proof outlined for the Picard's iterates is much more straightforward, having been proven mostly the same by most textbooks \cite{callier1991}. But regardless, it still stands that these proofs by themselves often present a significant barrier to understanding. Unlike theorems such as Poincaré-Bendixson, which are naturally visual, these rely heavily on abstract analysis that can be difficult to conceptualize. The overall intent of this manuscript was to add a little more intuition to the proofs already described in contemporary literature.

\end{document}